\documentclass[review]{elsarticle}

\usepackage{amsmath}

\usepackage{lineno,hyperref}
\modulolinenumbers[5]

\journal{Journal of \LaTeX\ Templates}

\newtheorem{theorem}{Theorem}

\numberwithin{equation}{section}

\begin{document}

\begin{frontmatter}

\title{ON GEOMETRY OF THE FIRST AND THE SECOND FUNDAMENTAL FORMS OF CANAL SURFACES}

\author{Y\i lmaz Tun\c{c}er}

\address{Usak University,
Science and Art Faculty, Mathematics Department, Usak TURKEY}

\ead{yilmaz.tuncer@usak.edu.tr}

\begin{abstract}
In this study, we analyze the general canal surfaces in terms of the
features that flat, II-flat, minimality and II-minimality under conditions K=0, H=0, K$_{II}$=0 and H$_{II}$ =0. Thus, we classified the
general non-degenerate canal surfaces according to their radiouses and the curvature of the centered curve.
\end{abstract}


\end{frontmatter}

\linenumbers

\section{Introduction}

For many years, surface theory has been a popular topic for many researchers
in many aspects. Besides the using curves and surfaces, canal surfaces are
the most popular in computer aided geometric design such that designing
models of internal and external organs, preparing of
terrain-infrastructures, constructing of blending surfaces, reconstructing
of shape, robotic path planning, etc. (see, \cite{Fr,Wm,Xu}). The knowledge
of first fundamental form I and second fundamental form II of a surface
facilitates the analysis and the classification of surface shape.
Especially, the geometry of the second fundamental form II has become an
important issue in terms of investigating intrinsic and extrinsic geometric
properties of the surfaces. Very recent results concerning the curvature
properties associated to II and other variational aspects can be found in
\cite{Sv}. One may associate to such a surface M geometrical objects
measured by means of its second fundamental form, as second mean curvature H$%
_{II}$ and second Gaussian curvature K$_{II}$, respectively. We are able to
compute K$_{II}$ and H$_{II}$\ of a surface by replacing the components of
the first fundamental form $E,F,G$ by the components of the second
fundamental form $e,f,g$ in Brioschi formula which is given by Francesco
Brioschi in the years of 1800's.

In this study, we investigated the flat, II-flat, minimal and II-minimal
canal surfaces in $IR^{3}$ by obtaining the curvatures K, H, K$_{II}$ and H$%
_{II}$.

Let M be a surface immersed in Euclidean 3-space, the first fundamental form
$I$ of the surface M is defined by
\begin{equation*}
I=Edu^{2}+2Fdudv+Gdv^{2}
\end{equation*}%
where $E=<M_{s},M_{s}>$, $F=<M_{s},M_{t}>$, $G=<M_{t},M_{t}>$. A surface is
called degenerate if it has the degenerate first fundamental form. The
second fundamental form $II$ of M is given by%
\begin{equation*}
II=edu^{2}+2fdudv+gdv^{2}
\end{equation*}%
where $e=\langle M_{ss},n\rangle $, $f=$ $\langle M_{st},n\rangle $, $%
g=\langle M_{tt},n\rangle $ and $n$ is the unit normal of M. The Gaussian
curvature K and the mean curvature H are given by
\begin{equation}
\text{K}=\frac{eg-f^{2}}{EG-F^{2}},  \label{01}
\end{equation}%
\begin{equation}
\text{H}=\frac{Eg-2Ff+Ge}{2(EG-F^{2})}  \label{02}
\end{equation}%
respectively. A regular surface is flat if and only if its Gaussian
curvature vanishes identically. A minimal surface in IR$^{3}$ is a regular
surface if its mean curvature vanishes identically\cite{Ga}.

Furthermore, the second Gaussian curvature K$_{II}$ of a surface is defined
by%
\begin{eqnarray}
\text{K}_{II} &=&\frac{1}{\left( \left\vert eg\right\vert -f^{2}\right) ^{2}}%
\left\{ p-q\right\} .  \notag \\
&&  \label{03}
\end{eqnarray}%
where
\begin{equation*}
p=\left\vert
\begin{array}{ccc}
-\frac{1}{2}e_{tt}+f_{st}-\frac{1}{2}g_{ss} & \frac{1}{2}e_{s} & f_{s}-\frac{%
1}{2}e_{t} \\
f_{t}-\frac{1}{2}g_{s} & e & f \\
\frac{1}{2}g_{t} & f & g%
\end{array}%
\right\vert
\end{equation*}%
and
\begin{equation*}
q=\left\vert
\begin{array}{ccc}
0 & \frac{1}{2}e_{t} & \frac{1}{2}g_{s} \\
\frac{1}{2}e_{t} & e & f \\
\frac{1}{2}g_{s} & f & g%
\end{array}%
\right\vert
\end{equation*}%
Similar to the variational characterization of the mean curvature H, the
curvature of the second fundamental form, denoted by H$_{II}$ is introduced
as a measure for the rate of change of the $II$ area under a normal
deformation. For details see \cite{Shv}. A surface is called II-flat if the
second Gaussian curvature vanishes identically\cite{Yoon}. Having in mind
the usual technique for computing the second mean curvature by using the
normal variation of the area functional one gets%
\begin{equation}
\text{H}_{II}=\text{H}-\frac{1}{2\sqrt{\left\vert \det II\right\vert }}%
\sum_{i,j}\frac{\partial }{\partial u^{i}}\left( \sqrt{\left\vert \det
II\right\vert }L^{ij}\frac{\partial }{\partial u^{j}}\left( \ln \sqrt{%
\left\vert K\right\vert }\right) \right)   \label{04}
\end{equation}%
where $u^{i}$ and $u^{j}$ stand for "$s$" and "$t$", respectively, and $%
\left( L^{ij}\right) =\left( L_{ij}\right) ^{-1},$ where $L_{ij}$ are the
coefficients of second fundamental forms\cite{Shv}. A surface is called
II-minimal if the second mean curvature vanishes identically\cite{Yoon}.

A general canal surface is an envelope of a 1-parameter family of surface.
The envelope of a 1-parameter family $s\longrightarrow S^{2}\left( s\right) $
of spheres in $IR^{3}$ is called a \textit{general} \textit{canal surface}%
\cite{Ga}.. The curve formed by the centers of the spheres \ is called
\textit{center curve }of the canal surface. The radius of general canal
surface is the function $r$ such that $r(s)$ is the radius of the sphere $%
S^{2}\left( s\right) .$ Suppose that the center curve of a canal surface is
a unit speed curve $\alpha :I\rightarrow IR^{3}$. Then the general canal
surface can be parametrized by the formula
\begin{equation}
C\left( s,t\right) =\alpha \left( s\right) -R\left( s\right) T-Q\left(
s\right) \cos \left( t\right) N+Q\left( s\right) \sin \left( t\right) B
\label{0}
\end{equation}%
where%
\begin{gather}
R\left( s\right) =r(s)r^{\prime }(s)  \label{1} \\
Q\left( s\right) =\pm r(s)\sqrt{1-r^{\prime }(s)^{2}}.  \label{2}
\end{gather}%
All the tubes and the surfaces of revolution are subclass of the general
canal surface.

\begin{theorem}
Let M be a canal surface. The center curve of M is a straight line if and
only if M is a surface of revolution for which no normal line to the surface
is parallel o the axis of revolution \cite{Ga}.
\end{theorem}

\begin{theorem}
The following conditions are equivalent for a canal surface M: \emph{i}. M
is a tube parametrized by (\ref{0}); \emph{ii}. the radius of M is constant;
\emph{iii}. the radius vector of each sphere in family that defines the
canal surface M meets the center curve orthogonally \cite{Ga}.
\end{theorem}

Coefficients of the first and the second fundamental forms of canal surface
are%
\begin{eqnarray}
E\left( s,t\right) &=&Q^{2}\kappa ^{2}\cos ^{2}\left( t\right) +p_{1}\kappa
\cos \left( t\right) +2QR\kappa \tau \sin (t)+p_{2}  \notag \\
F\left( s,t\right) &=&-Q\left( R\kappa \sin (t)+Q\tau \right)  \label{3a} \\
G\left( s,t\right) &=&Q^{2}  \notag
\end{eqnarray}%
and%
\begin{eqnarray}
e\left( s,t\right) &=&\frac{-1}{r\left( s\right) }\left\{ E-Q\kappa \cos
\left( t\right) -p_{5}\right\}  \notag \\
f\left( s,t\right) &=&\frac{-1}{r\left( s\right) }F\left( s,t\right)
\label{3b} \\
g\left( s,t\right) &=&\frac{-1}{r\left( s\right) }G\left( s,t\right)  \notag
\end{eqnarray}%
Let us take $\psi \left( s,t\right) =\det I$ and $\phi \left( s,t\right)
=\det II.$ Thus, we have
\begin{equation}
\phi \left( s,t\right) =\frac{1}{r^{2}}\left\{ \psi \left( s,t\right)
-Q^{3}\kappa \cos \left( t\right) -Q^{2}p_{5}\right\}  \label{3c}
\end{equation}%
\begin{eqnarray}
\psi \left( s,t\right) &=&Q^{2}\left\{ \kappa ^{2}\left( R^{2}+Q^{2}\right)
\cos ^{2}\left( t\right) +\kappa p_{1}\cos \left( t\right) +1-2R^{\prime
}+\left( R^{\prime }\right) ^{2}+\left( Q^{\prime }\right) ^{2}\right\} .
\notag \\
&&  \label{k2}
\end{eqnarray}%
and%
\begin{eqnarray}
p_{1} &=&2\left( Q-QR^{\prime }+Q^{\prime }R\right)  \notag \\
p_{2} &=&Q^{2}\tau ^{2}+R^{2}\kappa ^{2}+\left( R^{\prime }\right)
^{2}+\left( Q^{\prime }\right) ^{2}-2R^{\prime }+1  \notag \\
p_{3} &=&p_{1}-Q  \label{3d} \\
p_{4} &=&p_{2}-p_{5}  \notag \\
p_{5} &=&\left( R^{\prime }\right) ^{2}+\left( Q^{\prime }\right)
^{2}-2R^{\prime }+1+RR^{\prime \prime }+QQ^{\prime \prime }  \notag
\end{eqnarray}%
If $Q(s)=0$, then the first and the second fundamental forms are degenerate.
So the canal surface is degenerate surface and the radius is $r(s)=\pm s+c$.
Furthermore, in the case $\kappa (s)=0$ and $1-2R^{\prime }+\left( R^{\prime
}\right) ^{2}+\left( Q^{\prime }\right) ^{2}=0,$ the radius is
\begin{equation*}
r(s)=\sqrt{s^{2}-2c_{1}s+2c_{2}}.
\end{equation*}%
Let the center curve be $\alpha \left( s\right) =\left( s,0,0\right) $. Then
$T=e_{1}$, $N=e_{2}$ and $B=e_{3}$. Hence, $R\left( s\right) =s-c_{1}$ and $%
C\left( s,t\right) $ is the curve in the plane $x=c_{1}$. The conditions
that $r(s)\neq \pm s+c$ and $\left( \kappa (s)=0,r(s)\neq \sqrt{%
s^{2}-2c_{1}s+2c_{2}}\right) $ are the necessary conditions to define a
non-degenerate canal surface with the equation (\ref{1}). At this point, we
can write the following theorem.

\begin{theorem}
Let $M$ be a canal surface with the center curve $\alpha (s)$ and the radius
$r(s)$. If the center curve is a line then $M$ is a regular surface in IR$%
^{3}$ iff the radius is $r(s)\neq \pm s+c$ and $r(s)\neq \sqrt{%
s^{2}-2c_{1}s+2c_{2}}$.
\end{theorem}

Additionally, if $\phi \left( s,t\right) =0$ then $M$ has degenerate second
fundamental form. A canal surface has degenerate second fundamental form if
canal surface is a surface of revolution with the radious $r(s)=\sqrt{%
s^{2}-2c_{1}s+2c_{2}}$ or $r(s)=c_{1}s+c_{2}$.

\section{Specific Curvatures of General Canal Surfaces}

In this section, we obtained the curvatures K, H, K$_{II}$ and H$_{II}$ of
canal surface in IR$^{3}$ by using equations (\ref{01}), (\ref{02}), (\ref%
{03}), ( \ref{04})\ and we classify the canal surfaces under the conditions
for K$=0$, H$=0$, K$_{II}=0$ and H$_{II}=0$. From (\ref{01}), we obtaine the
Gauss curvature as a polynomial expression of cos$\left( t\right) $ such that%
\begin{eqnarray}
K(s,t) &=&\frac{-1}{\psi (s,t)r^{2}}\left\{ Q^{3}\kappa \cos \left( t\right)
+Q^{2}p_{5}-\psi (s,t)\right\}   \notag \\
&&  \label{k1}
\end{eqnarray}%
The condition for vanishing Gaussian curvature requires that $\kappa =0$ and
\begin{equation*}
RR^{\prime \prime }+QQ^{\prime \prime }=0.
\end{equation*}%
From (\ref{1}) and (\ref{2}) above equation turns to
\begin{equation}
r^{\prime \prime }(s)\left\{ r(s)r^{\prime \prime }(s)+\left( r^{\prime
}(s)\right) ^{2}-1\right\} =0.
\end{equation}%
The solutions are $r(s)=c\neq 0$ and $r(s)=c_{1}s+c_{2},$ $\left\vert
c_{1}\right\vert <1.$ By setting $\alpha \left( s\right) =(s,0,0),$ $T=e_{1}$%
, $N=e_{2}$ and $B=e_{3}$, (\ref{1}) is a cylinder
\begin{equation*}
C\left( s,t\right) =\left( s,\mp c\cos \left( t\right) ,\pm c\sin \left(
t\right) \right)
\end{equation*}%
and is a cone%
\begin{equation*}
C\left( s,t\right) =\left( \text{p}s-c_{1}c_{2},\mp \left(
c_{1}s+c_{2}\right) \sqrt{\text{p}}\cos \left( t\right) ,\mp \left(
c_{1}s+c_{2}\right) \sqrt{\text{p}}\sin \left( t\right) \right)
\end{equation*}%
respectively, where p$=1-\left( c_{1}\right) ^{2}$. Hence, we can state the
following theorem in the case $\kappa (s)=0$ and $Q(s)\neq 0$.

\begin{theorem}
Let M be a canal surface with the center curve $\alpha (s)$ and the radius $%
r(s)$ then M is a non-degenerate flat canal surface if M is either a
cylinder or a cone.
\end{theorem}

As in Gaussian curvature, from (\ref{02}), we obtaine the mean curvature in
a polynomial expression of $cos\left( t\right) $ as follow
\begin{eqnarray}
\text{H}(s,t) &=&\frac{1}{2\psi r^{2}}\left\{ Q^{3}\kappa \cos \left(
t\right) +Q^{2}p_{5}-2\psi \right\} .  \notag \\
&&  \label{h1}
\end{eqnarray}%
The condition for minimality requires that $\kappa =0$ and
\begin{equation}
\left( R^{\prime }\right) ^{2}+\left( Q^{\prime }\right) ^{2}-2R^{\prime
}-RR^{\prime \prime }-QQ^{\prime \prime }+1=0  \label{h2}
\end{equation}%
and from(\ref{1}) and (\ref{2}) and (\ref{h2}) turns to
\begin{equation*}
?\left( r(s)r^{\prime \prime }(s)+\left( r^{\prime }(s)\right) ^{2}-1\right)
\left( 2r(s)r^{\prime \prime }(s)+\left( r^{\prime }(s)\right) ^{2}-1\right)
=0.
\end{equation*}%
The only reel solution is $r(s)=\sqrt{s^{2}-2c_{1}s+2c_{2}}$. In this case, (%
\ref{1}) is degenerate surface. At this point, we can have the following
theorem.

\begin{theorem}
Let M be a canal surface with the center curve $\alpha (s)$ and the radius $%
r(s)$. Then there are no non-degenerate minimal canal surfaces in IR$^{3}.$
\end{theorem}

After long and calculations, we are able to obtaine the second Gauss and
second mean curvature by using (\ref{03})and (\ref{04}). In a polynomial
expression of cos$\left( t\right) $, the second Gauss curvature is
\begin{equation}
\text{K}_{II}(s,t)=\frac{-1}{4r^{5}\phi }\sum\limits_{i=0}^{2}n_{i}\cos
^{i}\left( t\right)   \label{kk1}
\end{equation}%
where the coefficients $n_{i}$ are
\begin{equation}
n_{2}=-r^{2}Q^{4}\kappa ^{2}  \label{kk2}
\end{equation}%
\begin{equation}
n_{1}=-2Q\left\{
\begin{array}{l}
r^{2}E_{tt}\kappa Q^{2}-2r^{2}F_{st}\kappa Q^{2}+2r^{2}\kappa Q^{3}Q^{\prime
\prime }-rr^{\prime \prime }\kappa Q^{4} \\
+rQF_{t}\left\{ r\kappa Q^{\prime }+Q(\kappa ^{\prime }r+\kappa r^{\prime
})\right\} -r^{2}\kappa Q^{2}\left( Q^{\prime }\right) ^{2} \\
-\frac{1}{2}rQ^{3}Q^{\prime }(2\kappa ^{\prime }r+\kappa r^{\prime
})+r^{2}\kappa (p_{5}Q^{2}-\psi ) \\
+\left( r^{\prime }\right) ^{2}\kappa Q^{4}+\frac{1}{2}rr^{\prime }\kappa
^{\prime }Q^{4}%
\end{array}%
\right\}   \label{kk3}
\end{equation}%
and%
\begin{eqnarray}
n_{0} &=&\left\{
\begin{array}{l}
2(Q^{2}p_{5}-\psi )\left\{ r^{\prime }Q(rQ^{\prime }-(r^{\prime
})+r^{2}(2F_{ts}-E_{tt}-2QQ^{\prime \prime })+rr^{\prime \prime
}Q^{2}\right\}  \\
+2rF_{t}\left\{ rFE_{t}+rQ^{2}E_{s}-2rFF_{s}+rQ(F\kappa \sin
(t)-(p_{5})^{\prime })\right\}  \\
-r^{2}Q^{2}(E_{t})^{2}+E_{t}\left\{ 2r^{2}QQ^{\prime }F-rr^{\prime
}Q^{2}F-2r^{2}Q^{3}\kappa \sin (t)\right\}  \\
+Q(r^{\prime }Q-2rQ^{\prime })\left\{
rQ^{2}E_{s}-2rFF_{s}-rQ(Q(p_{5})^{\prime }+F\kappa \sin (t))\right\}  \\
+2rr^{\prime }QQ^{\prime }F^{2}-4r^{2}\left( Q^{\prime }\right)
^{2}F^{2}-r^{2}Q^{4}\kappa ^{2}-2rr^{\prime }F_{t}Q^{2}(p_{5}-E)%
\end{array}%
\right\}   \notag \\
&&  \label{kk4}
\end{eqnarray}%
The condition for II-flatness requires that all coefficients $n_{i}$ are
zero in ( \ref{kk1}). From (\ref{kk2}) one has $\kappa =0,$ and \ in this
case, $n_{1}=0$. Thus, from equations (\ref{1}), (\ref{2}), (\ref{3a}) and
its derivatives with respect to s an t, $n_{0}=0$ turns to
\begin{gather}
0=rr^{\prime }r^{\prime \prime \prime }(r^{\prime }-1)^{2}(r^{\prime
}+1)^{2}-4r^{3}(r^{\prime \prime })^{4}-6r^{2}(r^{\prime \prime
})^{3}((r^{\prime })^{2}-1)  \label{kk7} \\
+2r(r^{\prime \prime })^{2}\left\{ 3(r^{\prime })^{2}-2(r^{\prime
})^{2}-1\right\} +(r^{\prime })^{2}r^{\prime \prime }\left\{ (r^{\prime
})^{4}-2(r^{\prime })^{2}+1\right\}   \notag
\end{gather}%
The reel solutions of (\ref{kk7}) are $r=c\neq 0$ and $r(s)=\pm s+c.$ Hence,
we can state the following theorem.

\begin{theorem}
Let $M$ be a non-degenere canal surface with the center curve $\alpha (s)$
and the radius $r(s)$. Then $M$ is a II-flat canal surface if the surface is
a cylinder.
\end{theorem}

The second mean curvature can be obtained by using (\ref{04}), in a
polynomial expression of cos$\left( t\right) $ with the aid of maple \textbf{%
prog. 2 }and we can give just coeffcients which we need. The second mean
curvature is
\begin{equation}
\text{H}_{II}(s,t)=\frac{1}{denom_{1}}\left\{ 4Q^{15}\kappa ^{5}\phi
^{2}\psi \cos ^{5}\left( t\right) +...+(...)\cos \left( t\right)
+...\right\} .  \label{hii1}
\end{equation}%
The condition for II-minimality requires $\kappa =0,$ in this case, all the
coefficients of the term cos$\left( t\right) $ are zero and the last term of
(\ref{hii1}) is too long. From (\ref{1}), (\ref{2}) and prog.1 and prog.2,
we are able to the differential equation obtained from the last term of (\ref%
{hii1}). The solutions are $r(s)=c$, $r(s)=c_{1}s+c_{2}$ and $r(s)=\pm \sqrt{%
s^{2}-2c_{1}s+2c_{2}}.$Thus we can give the following theorem.

\begin{theorem}
Let M be a non-degenere canal surface with the center curve $\alpha (s)$ and
the radius $r(s)$. Then M is a II-minimal canal surface if M is either a
tubular surface with the centered line or a cone.
\end{theorem}

On the other hand, we can consider (\ref{0}) as
\begin{equation}
\text{H}_{II}=\text{H}-\frac{\left( \delta _{1}+\delta _{2}\right) }{2\sqrt{%
\mu _{1}\phi }}  \label{3}
\end{equation}%
where%
\begin{eqnarray*}
\delta _{1} &=&\frac{\partial }{\partial s}\left( \sqrt{\mu _{1}\phi }L^{11}%
\frac{\partial }{\partial s}\left( \ln \sqrt{\mu _{2}K}\right) +\sqrt{\mu
_{1}\phi }L^{12}\frac{\partial }{\partial t}\left( \ln \sqrt{\mu _{2}K}%
\right) \right)  \\
\delta _{2} &=&\frac{\partial }{\partial t}\left( \sqrt{\mu _{1}\phi }L^{21}%
\frac{\partial }{\partial s}\left( \ln \sqrt{\mu _{2}K}\right) +\sqrt{\mu
_{1}\phi }L^{22}\frac{\partial }{\partial t}\left( \ln \sqrt{\mu _{2}K}%
\right) \right)
\end{eqnarray*}%
and
\begin{equation*}
\mu _{1}=\left\{
\begin{array}{c}
1 \\
-1%
\end{array}%
\begin{array}{c}
;\text{if }\phi >0 \\
;\text{if }\phi <0%
\end{array}%
\right. \text{ \ and }\mu _{2}=\left\{
\begin{array}{c}
1 \\
-1%
\end{array}%
\begin{array}{c}
;\text{if }K>0 \\
;\text{if }K<0%
\end{array}%
\right. .
\end{equation*}%
If $\delta _{1}+\delta _{2}=0$ then, M has the same mean curvature value
according to the first and second fundamental form. With the aid of prog.3,
we obtain $\delta _{1}+\delta _{2}$ as follow.%
\begin{equation*}
\delta _{1}+\delta _{2}=\frac{1}{denom_{2}}\sum\limits_{i=0}^{3}z_{i}\cos
^{i}\left( t\right)
\end{equation*}%
In the case $\delta _{1}+\delta _{2}=0,$ all the coefficients $z_{i}$ are
zero. For $z_{3}=0,$ we have
\begin{equation*}
-\mu _{1}Q^{7}r^{2}\kappa ^{3}\left\{ \psi ^{2}+2\psi \psi _{tt}-2(\psi
_{t})^{2}\right\} =0.
\end{equation*}%
If $\kappa =0$ and $\psi ^{2}+2\psi \psi _{tt}-2(\psi _{t})^{2}\neq 0$ then $%
z_{3}=z_{2}=z_{1}=z_{0}=0.$ If $\kappa \neq 0$ and $\psi ^{2}+2\psi \psi
_{tt}-2(\psi _{t})^{2}=0$ then, from (\ref{3a}), (\ref{3b}), (\ref{3c}), (%
\ref{k2}) and with the aid of prog.3, we have the following%
\begin{equation}
r^{2}Q^{4}\kappa ^{4}\cos ^{4}\left( t\right) +(...)\cos ^{2}\left( t\right)
+(...)\cos \left( t\right) +...=0.  \label{33}
\end{equation}%
It is obvious that (\ref{3a}) don't satisfy. If $\kappa =0$ ve $\psi
^{2}+2\psi \psi _{tt}-2(\psi _{t})^{2}=0$ then, with the aid of prog.1 and
prog.3, we get%
\begin{equation*}
Q^{4}\left( p_{2}\right) ^{2}=0
\end{equation*}%
and by using $p_{2}$ in (\ref{3d}) and from (\ref{1}) and (\ref{2}) we get
the solutions $r=\pm s+c_{1}$\ and $r=\pm \sqrt{s^{2}-2c_{1}s+2c_{2}}.$Thus,
we have the following theorem.

\begin{theorem}
\bigskip Let M be a non-degenere canal surface with the center curve $\alpha
(s)$ and the radius $r(s)$. Then, the first and the second mean curvatures
of M are equal if M is surface of revolution.
\end{theorem}

Similarly, we can consider the canal surfaces whose the Gaussian curvatures
according to the first and second fundamental forms are equal. From (\ref{k1}%
) and (\ref{kk1}), for K$(s,t)-$K$_{II}(s,t)=0$%
\begin{equation}
\psi n_{2}\cos ^{2}\left( t\right) +\left( \psi n_{1}-4\phi r^{3}Q^{3}\kappa
\right) \cos \left( t\right) +\psi n_{0}+4\phi r^{3}(\psi -Q^{2}p_{5})=0
\label{34}
\end{equation}%
where $n_{0},n_{1}$ and $n_{2}$ are as in (\ref{kk2}),(\ref{kk3}) and (\ref%
{kk4}). The equality (\ref{34}) requires $\kappa =0$. For $\kappa =0,$
\begin{equation*}
\psi n_{0}+4\phi r^{3}(\psi -Q^{2}p_{5})=0
\end{equation*}%
turns to a too long differential equation. With the aid of prog.4, we
obtaine the solutions that $r=\pm s+c_{1}$\ and $r=\pm \sqrt{%
s^{2}-2c_{1}s+2c_{2}}.$Thus, we have the following theorem.

\begin{theorem}
\bigskip Let M be a non-degenere canal surface with the center curve $\alpha
(s)$ and the radius $r(s)$. Then, there is no canal surface whose the first
and the second Gaussian curvatures of M are equal.
\end{theorem}

\textbf{Prog.1}

p1(s):=2*Q(s)+2*R(s)*diff(Q(s),s)-2*Q(s)*diff(R(s),s):

p2(s):=(Q(s)\symbol{94}2)*(tau(s)\symbol{94}2)+(R(s)\symbol{94}2)*(kappa(s)%
\symbol{94}2)+diff(R(s),s)\symbol{94}2

+diff(Q(s),s)\symbol{94}2-2*diff(R(s),s)+1:

p3(s):=p1(s)-Q(s):\qquad

p5(s):=diff(R(s),s)\symbol{94}2+diff(Q(s),s)\symbol{94}%
2-2*diff(R(s),s)+1+R(s)*diff(diff(R(s),s),s)

+Q(s)*diff(diff(Q(s),s),s):

p4(s):=p2(s)-p5(s):

E(s,t):=(Q(s)\symbol{94}2)*(kappa(s)\symbol{94}2)*(cos(t))\symbol{94}%
2+p1(s)*kappa(s)*cos(t)

+2*Q(s)*R(s)*kappa(s)*tau(s)*sin(t)+p2(s):

F(s,t):=-Q(s)*R(s)*kappa(s)*sin(t)-G(s)*tau(s):

\textbf{Prog.2}

e:=(-1/r(s))*(E(s,t)-Q(s)*kappa(s)*cos(t)-p5(s)):

f:=(-1/r(s))*F(s,t):

g:=(-1/r(s))*G(s):

L11:=g/(e*g-f\symbol{94}2):

L12:=-f/(e*g-f\symbol{94}2):

L21:=-f/(e*g-f\symbol{94}2):

L22:=e/(e*g-f\symbol{94}2):

G(s):=(Q(s))\symbol{94}2;

K:=1/psi(s,t)*(-Q(s)\symbol{94}3*kappa(s)*cos(t) -Q(s)\symbol{94}%
2*p5(s)+psi(s,t))/r(s)\symbol{94}2:

H:=-1/2*(-Q(s)\symbol{94}3*kappa(s)*cos(t) -Q(s)\symbol{94}2*p5(s)+
2*psi(s,t)\symbol{94}2)/psi(s,t)/r(s):

delta1:=simplify(diff(((mu1*(phi(s,t)))\symbol{94}(1/2))*L11*diff(ln((mu2*K)%
\symbol{94}(1/2)),s)

+((mu1*(phi(s,t)))\symbol{94}(1/2))*L12*diff(ln((mu2*K)\symbol{94}%
(1/2)),t),s)):

delta2:=simplify(diff(((mu1*(phi(s,t)))\symbol{94}(1/2))*L21*diff(ln((mu2*K)%
\symbol{94}(1/2)),s)

+((mu1*(phi(s,t)))\symbol{94}(1/2))*L22*diff(ln((mu2*K)\symbol{94}%
(1/2)),t),t)):

H2:=simplify(H-(1/(2*(mu1*(phi(s,t)))\symbol{94}(1/2)))*(delta1+delta2)):

Hoo:=subs(cos(t)=A,H2):

Ho:=subs(sin(t)=B,Hoo):

simplify(coeff(numer(Ho),A,5),'size');

\textbf{Prog.3}

e:=(-1/r(s))*(E(s,t)-Q(s)*kappa(s)*cos(t)-p5(s)):

f:=(-1/r(s))*F(s,t):

g:=(-1/r(s))*G(s):

L11:=g/(phi(s,t)):

L12:=-f/(phi(s,t)):

L21:=-f/(phi(s,t)):

L22:=e/(phi(s,t)):

G(s):=(Q(s))\symbol{94}2;

phi(s,t):=e*g-f\symbol{94}2:

K:=simplify(((e*g-f\symbol{94}2)/(psi(s,t)))):

delta1:=simplify(diff(((mu1*(phi(s,t)))\symbol{94}(1/2))*L11*diff(ln((mu2*K)%
\symbol{94}(1/2)),s)

+((mu1*(phi(s,t)))\symbol{94}(1/2))*L12*diff(ln((mu2*K)\symbol{94}%
(1/2)),t),s)):

delta2:=simplify(diff(((mu1*(phi(s,t)))\symbol{94}(1/2))*L21*diff(ln((mu2*K)%
\symbol{94}(1/2)),s)

+((mu1*(phi(s,t)))\symbol{94}(1/2))*L22*diff(ln((mu2*K)\symbol{94}%
(1/2)),t),t)):

z:=simplify(delta1+delta2):

zoo:=subs(cos(t)=A,z):

zo:=subs(sin(t)=B,zoo):

simplify(coeff(numer(zo),A,0),'size');

\textbf{Prog.4}

kappa(s):=0:

tau(s):=0:

R(s):=r(s)*diff(r(s),s):

Q(s):=r(s)*((1-diff(r(s),s)\symbol{94}2)\symbol{94}(1/2)):

p2(s):=(Q(s)\symbol{94}2)*(tau(s)\symbol{94}2)+(R(s)\symbol{94}2)*(kappa(s)%
\symbol{94}2)+diff(R(s),s)\symbol{94}2

+diff(Q(s),s)\symbol{94}2-2*diff(R(s),s)+1:

p5(s):=diff(R(s),s)\symbol{94}2+diff(Q(s),s)\symbol{94}%
2-2*diff(R(s),s)+1+R(s)*diff(diff(R(s),s),s)

+Q(s)*diff(diff(Q(s),s),s):

E(s,t):=(Q(s)\symbol{94}2)*(kappa(s)\symbol{94}2)*(cos(t))\symbol{94}%
2+p1(s)*kappa(s)*cos(t)

+2*Q(s)*R(s)*kappa(s)*tau(s)*sin(t)+p2(s):

F(s,t):=-Q(s)*R(s)*kappa(s)*sin(t)-G(s)*tau(s):

G(s):= Q(s)\symbol{94}2:

psi(s,t):=E(s,t)*G(s)-F(s,t)\symbol{94}2:

e:=(-1/r(s))*(E(s,t)-Q(s)*kappa(s)*cos(t)-p5(s)):

f:=(-1/r(s))*F(s,t):

g:=(-1/r(s))*G(s):

phi:=e*g-f\symbol{94}2:

es:=diff(e,s):

fs:=diff(f,s):

gs:=diff(g,s):

ess:=diff(diff(e,s),s):

fss:=diff(diff(f,s),s):

gss:=diff(diff(g,s),s):

et:=diff(e,t):

ft:=diff(f,t):

gt:=diff(g,t):

ett:=diff(diff(e,t),t):

ftt:=diff(diff(f,t),t):

gtt:=diff(diff(g,t),t):

est:=diff(diff(e,s),t):

fst:=diff(diff(f,s),t):

gst:=diff(diff(g,s),t):

ets:=diff(diff(e,t),s):

fts:=diff(diff(f,t),s):

gts:=diff(diff(g,t),s):

V1 :=
Matrix([[(-ett/2)+fst-(gss/2),(es/2),fs-(et/2)],[ft-(gs/2),e,f],[gt/2,f,g]]):

V2 := Matrix([[0,et/2,gs/2],[et/2,e,f],[gs/2,f,g]]):

v1 := LinearAlgebra:-Determinant(V1):

v2 := LinearAlgebra:-Determinant(V2):

K2:=simplify((v1-v2)/(phi)):

K:=1/psi(s,t)*(-Q(s)\symbol{94}3*kappa(s)*cos(t) -Q(s)\symbol{94}%
2*p5(s)+psi(s,t))/r(s)\symbol{94}2:

simplify(coeff(numer(subs(sin(t)=B,subs(cos(t)=A, K- K2))),A,0),'size');

\end{document}